\documentclass[11pt]{article}
\usepackage{CJK}
\usepackage{amsmath,amsfonts,mathrsfs,amssymb}
\usepackage{indentfirst}
\numberwithin{equation}{section}
\usepackage{color}

\usepackage[pagewise]{lineno}

\numberwithin{equation}{section}

\newcommand{\HT}{\CJKfamily{hei}}

\def\lb{\label}

\def\be{\begin{equation}}
\def\ee{\end{equation}}
\def\bea{\begin{eqnarray}}
\def\eea{\end{eqnarray}}
\def\bes{\begin{eqnarray*}}
\def\ees{\end{eqnarray*}}

\def\y{\begin{eqnarray*}}
\def\ey{\end{eqnarray*}}

\begin{document}
\bibliographystyle{amsplain}
\title{\HT {The ground state solution for  Kirchhoff-Schr\"odinger equations with singular
exponential nonlinearities in $\mathbb{R}^{4}$  }}
\author{\small {Yanjun Liu$^{a, ^{\ast}}$ ~~~~Shijie Qi$^{b}$ }\\
\small $^{a}$School of Mathematical Sciences, Nankai University, Tianjin 300071, P. R. China\\
\small $^{b}$Department of Mathematical Science, Tsinghua University, Beijing  100084, P. R. China}
\date{}
\maketitle \footnote[0] {$^{\ast}$ Corresponding author.}
\maketitle \footnote[0] {Email: liuyj@mail.nankai.edu.cn(Y. Liu); qishj2019@tsinghua.edu.cn(S. Qi).}

\noindent{\small
{\bf Abstract:} In this paper, we consider the following  singular  Kirchhoff-Schr\"odinger problem
$$M\bigg(\int_{\mathbb{R}^{4}}|\Delta u|^2+V(x)u^2 dx\bigg)(\Delta^{2}u+V(x)u) =\frac{f(x, u)}{|x|^{\eta}}~~\operatorname{in}~~ \mathbb{R}^{4}, \eqno(P_{\eta})$$
where $0<\eta<4$, $M$ is a Kirchoff-type function and $V(x)$ is a continuous function with positive lower bound, $f(x, t)$  has an critical exponential  growth behavior at infinity. Using  singular  Adams  inequality and variational techniques,  we get the existence of ground state solutions for $(P_{\eta})$. Moreover, we also get the same results without the Ambrosetti-Rabinowitz (AR) condition.


   }

\section{Introduction and main results}
 The relevant problems involving powers of the Laplacian started with \cite{111}-\cite{222}. In conformal geometry, there has been considerable interest in the Paneitz operator which enjoys the property of conformal invariance. In
$\mathbb{R}^{4}$, the Paneitz operator is the biharmonic operator $\Delta^2$, this can be refered to \cite{F. Sani}. Recently, Zhang and Chen in \cite{L. Chen} established  a sharp concentration-compactness principle associated with the singular Adams inequality on the second-order Sobolev spaces in $\mathbb{R}^{4}$, and  moreover,they consider the following problem:
\be \Delta^{2}u+V(x)u =\frac{f(x, u)}{|x|^{\eta}}~~\operatorname{in}~~ \mathbb{R}^{4}, \label{1} \ee
where $V(x)$ has a positive lower bound and $0<\eta <4$, they got a ground state solution of \eqref{1} under the A-R  condition. In \cite{GS}, let $\Omega\subset \mathbb{R}^{N}$ is a bounded domain with smooth boundary, $m$ is a integer and $N \geq 2m \geq 2$, the authors considered the following Kirchhoff problem
 \begin{equation}
 \left\{  \begin{array}{l}
         -M\bigg(\int_{\Omega}|\nabla^m u|^{\frac{N}{m}}dx\bigg)\Delta^m_{\frac{N}{m}}u =\frac{f(x, u)}{|x|^{\eta}}~~~~\operatorname{in}~~\Omega, \\
          u=\nabla u= \nabla^2 u=\cdots  \nabla^{m-1} u=0 ~~\operatorname{on}~~\partial \Omega, \\
         \end{array}
         \right.\label{2}
 \end{equation}
where $0\leq \eta <N$, $M$ is a Kirchoff-type function and $b(x)$ is a continuous function with positive lower bound, $f(x, t)$  has an critical exponential  growth behavior at infinity.

Since we will work with exponential critical growth, we need to review the Trudinger--Moser inequality and Adams inequality, the latter is a generalization of the former and more details are as follows: On one hand, let $\Omega$ denotes a smooth bounded domain in $\mathbb{R}^{N}(N\geq 2)$, N. Trudinger \cite{N.S. Trudinger} proved that there exists $\alpha>0$ such that $W^{1, N}_{0}(\Omega)$ is embedded in the Orlicz space $L_{\varphi_{\alpha}}(\Omega)$ determined by the Young function $\varphi_{\alpha}(t)=e^{\alpha|t|^{\frac{N}{N-1}}}$, it was sharpened by J. Moser \cite{J. Moser} who found the best exponent $\alpha$. On the other hand,the Trudinger--Moser inequality was extended for unbounded domains by D. M. Cao \cite{D.M. Cao} in $\mathbb{R}^{2}$ and for any dimension $N\geq2$ by J. M. do $\acute{O}$ \cite{J. M. do1}. Moreover, J. M. do $\acute{O}$ et al. \cite{J. M. do2} established a sharp concentration-compactness principle associated with the singular Trudinger--Moser inequality in $\mathbb{R}^{N}$. For more results concerning the Trudinger-Moser inequality and its application in $N$-Laplacian equations, one can refer to \cite{Adimurthi2,Yang,Lu,Zhang1,Li} and the references therein.

For Adams type inequality, let $\Omega \subset \mathbb{R}^{4}$ be a smooth bounded domain. D. Adams \cite{D. Adams} derives
$$\sup_{u\in W_{0}^{2, 2}(\Omega), \int_{\Omega}|\Delta u|^{2}dx\leq 1}\int_{\Omega}e^{32\pi^{2}u^{2}}dx<\infty,$$
which was extended by C. Tarsi \cite{C. Tarsi}, i.e.
$$\sup_{u\in W^{2, 2}(\Omega)\bigcap W_{0}^{1, 2}(\Omega), \int_{\Omega}|\Delta u|^{2}dx\leq 1}\int_{\Omega}e^{32\pi^{2}u^{2}}dx<\infty,$$
B. Ruf and F. Sani \cite{B. Ruf} extended the Adams inequality to $\mathbb{R}^{4}$, namely
$$\sup_{u\in W^{2, 2}(\mathbb{R}^{4}), \int_{\mathbb{R}^{4}}(-\Delta u+u)^{2}dx\leq 1}\int_{\mathbb{R}^{4}}(e^{32\pi^{2}u^{2}}-1)dx<\infty,$$
 where $32\pi^{2}$ is the best constant. In order to apply this inequality to partial differential equation more reasonably, Y. Yang in \cite{Y. Yang} proves the following singular Adams inequality:\\
 \textbf{Theorem A.} Suppose  $0\leq \eta <4$, $\tau, \sigma$ are two positive constants. Then
\be \sup_{u\in W^{2, 2}(\mathbb{R}^{4}), \int_{\mathbb{R}^{4}}(|\Delta u|^{2}+\tau |\nabla u|^{2}+\sigma u^{2})dx\leq 1}\int_{\mathbb{R}^{4}}\frac{e^{\alpha u^{2}}-1}{|x|^{\eta}}dx<\infty,   \label{3}\ee
 where $\alpha\leq32\pi^{2}(1-\frac{\eta}{4})$ is the sharp constant. If $\alpha>32\pi^{2}(1-\frac{\eta}{4})$, then the supremum is infinite.

In \cite{Z. Yang}, Li and Yang  studied the following Schr$\ddot{o}$dinger-Kirchhoff type equation
\begin{equation}
  \left\{  \begin{array}{l}
        \bigg(\int\limits_{\mathbb{R}^{N}}(|\nabla u|^N+V(x)|u|^N)dx\bigg)^k(-\Delta_N u+V(x)|u|^{N-2}u)=\lambda A(x)|u|^{p-2}u+f(u), \\
        u\in W^{1, N}(\mathbb{R}^{N}),
   \end{array}
   \right. \lb{(1.3)}
 \end{equation}
 where $\Delta_N u=\texttt{div}(|\nabla u|^{N-2}\nabla u)$, $k > 0$, $V : \mathbb{R}^{N} \rightarrow (0, \infty)$ is continuous,
$\lambda > 0$ is a real parameter, $A$ is a positive function in $L^{\frac{p}{p-q}}(\mathbb{R}^{N})$ and $f$ satisfies exponential
growth. They derived two nontrivial solutions of \eqref{(1.3)} as the parameter $\lambda $ small enough. Indeed,
 suppose $\Omega \subset  \mathbb{R}^{N}$ is
a bounded smooth domain, the above problems is related to the stationary analogue of the equation
\be  u_{tt}-\bigg(a+b\int_\Omega |\nabla u|^2dx\bigg)\Delta u=f(x, u)  \lb{(1.4)}\ee
proposed by Kirchhoff in \cite{Kirchhoff} as an extension of the classical D'Alembert's wave equation for free vibrations
of elastic strings. In \cite{J.L. Lions}, Lions proposed an abstract framework for the problem and after that, problem
\eqref{(1.3)} began to receive a lot of attention. In \cite{FZ}, the authors  studied the following Schr$\ddot{o}$dinger-Kirchhoff type equation
\be M\bigg(\int_{\mathbb{R}^{2}}|\nabla u|^2+V(x)u^2 dx\bigg)(-\Delta u+V(x)u) =A(x)f(u)~~\operatorname{in}~~ \mathbb{R}^{2} \ee
where  $M$ is a Kirchoff-type function and $V(x)\geq V_0$ is a continuous function, $A$ is
locally bounded and the function $f$ has critical exponential growth. Applying variational methods beside a new Trudinger-Moser type inequality, they get the existence of ground state solution.  Moreover,  in the the local case
$M \equiv 1$, they also get some relevant results.

In this paper, we consider the following  singular biharmonic Kirchhoff-Schr$\ddot{o}$dinger problem
$$M\bigg(\int_{\mathbb{R}^{4}}|\Delta u|^2+V(x)u^2 dx\bigg)(\Delta^{2}u+V(x)u) =\frac{f(x, u)}{|x|^{\eta}}~~\operatorname{in}~~ \mathbb{R}^{4}, \eqno(P_{\eta})$$
where $0<\eta<4$, $M$ is a Kirchoff-type function and $V(x)\geq V_0$ is a continuous function, $f(x, t)$  has an critical exponential  growth behavior at infinity. Using  singular  Adams  inequality and variational techniques,  we get the existence of ground state solutions for $(P_{\eta})$.

Let $\mathfrak{M}(t)=\int_0^{t}M(s)ds$, we assume that $M : \mathbb{R}^{+} \rightarrow \mathbb{R}^{+}$ is a continuous function with $M(0) = 0$, and satisfies

$(M_1)$  $M_0=\inf_{t\geq 0}M(t)>0$;

$(M_2)$  for any $t_1$, $t_2 \geq 0$, it holds
                      $$\mathfrak{M}(t_1+t_2)\geq \mathfrak{M}(t_1)+\mathfrak{M}(t_2);$$

$(M_3)$  $\frac{M(t)}{t}$ is decreasing in $(0, \infty)$;\\
\textbf{Remark 1.1.}  By $(M_3)$, we can obtain that
    $2\mathfrak{M}(t) - M(t)t$ is nondecreasing for $t > 0$,
In particular,
     \be 2\mathfrak{M}(t) - M(t)t\geq 0~~~  \forall t \geq 0.  \label{1.4}\ee

 we require that $f(x, t)=0$ for all $(x, t)\in \mathbb{R}^{4}\times (-\infty, 0]$. Furthermore, we assume the function $f$ satisfying:

$(f_{1})$ $f$ is a continuous function and  $f(x, t)>0$ for all $t> 0$.

$(f_{2})$  There exist constants $\alpha_{0}, c_{1}, c_{2} > 0$ such that for all\ $(x, t)\in \mathbb{R}^{4}\times \mathbb{R}^{+}$,
        $$f(x, t)\leq c_{1}|t|^{3}+c_{2}(e^{\alpha_{0} t^2-1}),$$

$(f_{3})$ There exist $R_0>0$ and $\mu >4$, for any $x\in \Omega$,
     \be \mu F(x, t)\leq t f(x, t)~~~~~\forall |t|\geq R_0.   \lb{1.11} \ee
where $F(x, t)=\int_0^\infty f(x, s) ds$. This is the so-called the Ambrosetti-Rabinowitz (AR) condition.\\
We also give the following conditions on the potential $V(x)$ :

$(V_{1})$  $V$ is a continuous function satisfying $V(x) \geq V_{0 }> 0$;

Define a function space
     $$E=\{u \in  W^{2, 2}(\mathbb{R}^{4}):  \int_{\mathbb{R}^{4}}(|\Delta u|^{2}+V(x)|\nabla u|^{2})dx<\infty\},$$
which be equipped with the norm
$$\|u\|_{E}=\left(\int_{\mathbb{R}^{4}}(|\Delta u|^{2}+V(x)|\nabla u|^{2})dx\right)^{\frac{1}{2}}, $$
then the assumption\ $(V_1)$ implies\ $E$ is a reflexive \ Banach space. For any $p\ge 2$, we define
  $$S_{p}=\inf_{u\in E\backslash\{0\}}\frac{\| u \|_{E}}{(\int_{\mathbb{R}^{4}}\frac { |u|^{p}}{|x|^{\eta}}dx)^{\frac{1}{p}}}$$ and

      $$\lambda_{\eta}=\inf_{u\in E\backslash\{0\}}\frac{\| u \|_{E}^{4}}{\int_{\mathbb{R}^{4}}\frac { |u|^{4}}{|x|^{\eta}}dx}.$$
The continuous embedding of\ $E\hookrightarrow W^{2, 2}(\mathbb{R}^{4})\hookrightarrow L^{2 p}(\mathbb{R}^{4})(p \geq 2)$ and\ H\"{o}lder inequality implies
 \begin{align*}
    &\int_{\mathbb{R}^{4}}\frac{|u|^{2p}}{|x|^{\eta}}dx\\
    \leq& \int_{\{|x|>1\}}|u|^{2p}dx
    +\left(\int_{\{|x|\leq 1\}}|u|^{2p t'}dx\right)^{\frac{1}{t'}}\left(\int_{\{|x|\leq 1\}}\frac{1}{|x|^{\eta t}}dx\right)^{\frac{1}{t}}\\
    \leq& C \|u\|_{E}^{2p}.
   \end{align*}
where\ $1/t + 1/t' = 1$ and $t>1$ such that $\eta t<4$. Thus we have $S_p > 0$. We now  introduce the following two conditions.

$(f_{4})$ $\limsup_{t\rightarrow 0^{+}}\frac{2F(x, t)}{t^{3}}<\mathfrak{M}(1)\lambda_{\eta}~~~~\operatorname{uniformly~in}~\mathbb{R}^{4}$.

$(f_{5})$   $\frac{f(x, t)}{t^3}$ is  increasing in $t > 0$.

Our main results can be stated as follows:\\
\textbf{Theorem 1.1.}  Suppose $V$ satisfies $(V_{1})$, $f$ satisfies $(f_{1})-(f_{5})$. Furthermore we
assume

$(f_{6})$ There exist constants $p > 4$ and $C_{p}$ such that for all $(x, t)\in \mathbb{R}^{4} \times [0, \infty)$
     $$f(x, t)\geq C_{p}t^{p-1},$$
where
$$C_{p}:=\inf\left\{C>0: p \mathfrak{M}(t^{2}S_{p}^{2})-2Ct^{p}< p\mathfrak{M}\big(((1-\frac{\eta}{4})\frac{32 \pi^2 }{\alpha_{0}})\big)\right\}.$$
Then the problem $(P_{\eta})$ has a nontrivial nonnegative ground state solution in $E$.

Now  instead the condition  $(f_3)$,    we assume that

$(f'_{3})$ $\lim_{|t|\rightarrow +\infty}\frac{F(x, t)}{|t|^{4}}=\infty$ uniformly on $x \in \mathbb{R}^{4}$. We derive the the results without the Ambrosetti-Rabinowitz (AR) condition.\\
\textbf{Theorem 1.2.}  Suppose $V$ satisfies $(V_{1})$, $f$ satisfies $(f_{1})-(f_{2})$, $(f'_{3})$ and ($f_{4})-(f_{6})$.
Then the  problem  $(P_{\eta})$ possesses a positive ground state solution.\\

This paper is organized as follows: In Section 2,  we give some
preliminary results. In Section 3, we study the functionals and compactness analysis. In section 4, we prove Theorem 1.1. In section 5, we study the results without the Ambrosetti-Rabinowitz (AR) condition.
\section{ Preliminaries}
In this section we will give some preliminaries for our use later.\\
\textbf{Lemma 2.1 (see \cite{L. Chen})} Suppose $q\geq 2$ and $0<s<4$. Then $E$ can be compactly embedded
into $L^{q}(\mathbb{R}^{4}, |x|^{-s}dx)$.\\
\textbf{Lemma 2.2.} Let $\beta>0, 0<\eta <4$ and $\|u\|_{E}\leq T$ such that $\beta T ^{2}<32\pi^{2}(1-\frac{\eta}{4})$ and $q>2$, then
 $$\int_{\mathbb{R}^{4}}\frac{e^{\beta|u|^{2}}-1}{|x|^{\eta}}|u|^{q}dx\leq C(\beta)\|u\|_{E}^{q}.$$
\textbf{Proof.} Set $R(\beta, u)=e^{\beta u^{2}}-1$, using the H\"{o}der inequality, we have
\begin{equation*}
\left.
\begin{aligned}[b]
\int_{\mathbb{R}^{4}}\frac{R(\beta, u)}{|x|^{\eta}}|u|^{q}dx
    \leq &\int_{|u|\leq 1}\frac{R(\beta, u)}{|x|^{\eta}}|u|^{q}dx +\int_{|u|> 1}\frac{R(\beta, u)}{|x|^{\eta}}|u|^{q}dx \\
          \leq&R(\beta, 1)\int_{|u|\leq 1}\frac{|u|^{q}}{|x|^{\eta}}dx+ \bigg(\int_{\mathbb{R}^{4}}\frac{R(p\beta, u)}{|x|^{\eta}}dx\bigg)^{\frac{1}{p}}\bigg(\int_{\mathbb{R}^{4}}\frac{|u|^{qp'}}{|x|^{\eta}}dx\bigg)^\frac{1}{p'} \\
         \leq&R(\beta, 1)\|u\|_E^q+ \bigg(\int_{\mathbb{R}^{4}}\frac{R(p\beta, u)}{|x|^{\eta}}dx\bigg)^{\frac{1}{p}}\|u\|_E^q
\end{aligned}
\right.
\end{equation*}
where the last inequality is a direct consequence of Lemma 2.1. Choosing $p>1$ is sufficiently close $1$ such that $\beta p T^{2}<32\pi^{2}(1-\frac{\eta}{4})$, $\frac{1}{p}+\frac{1}{p'}=1$, Then the result can be derived from Theorem A.   ~~~~~~~~$\hfill\Box$\\
\textbf{Lemma 2.3.}  If $(f_5)$ holds, then for all $x\in \mathbb{R}^4$, we have that
           $H (x, t) = tf (x, t)-4F (x, t)$
is  increasing in $t > 0$.\\
\textbf{Proof.} Let $0 < t_1 < t_2$ be fixed. It follows from $(f_5)$ that
      \be t_1f (x, t_1)-4F (x, t_1)<\frac{f(x, t_2)}{t_2^3}t_1^4-4F (x, t_2)+4\int_{t_1}^{t_2}f(x, s)ds.  \label{2.1}\ee
On the other hand,
  \be   4\int_{t_1}^{t_2}f(x, s)ds<4\frac{f(x, t_2)}{t_2^3}\int_{t_1}^{t_2}s^3 ds=\frac{f(x, t_2)}{t_2^3}(t_2^4-t_1^4).   \label{2.2}\ee
From \eqref{2.1} and \eqref{2.2}, we derive that
      $$t_1f (x, t_1)-4F (x, t_1)<  t_2f (x, t_2)-4F (x, t_2).$$
This completes the proof. $\hfill\Box$
\section{Mountain pass geometry and minimax estimates}
We say that\ $u\in E$ is a weak solution of problem\ $(P_{\eta})$ if  for all\ $\phi \in E$,
$$M(\|u\|_{E}^{2})\int_{\mathbb{R}^{4}}(\triangle u \triangle \phi+V(x)u \phi )dx-\int_{\mathbb{R}^{4}}\frac{f(x, u)}{|x|^{\eta}}\phi dx=0.$$

Define the functional $I: E\rightarrow \mathbb{R}$ by
$$I(u)=\frac{1}{2}\mathfrak{M}(\|u\|_{E}^{2}) -\int_{\mathbb{R}^{4}}\frac{F(x, u)}{|x|^{\eta}}dx.  \eqno(3.1)$$
where\ $F(x, t)=\int_{0}^{t}f(x, s)ds$. $I$ is well defined and \ $I\in C^{1}(E, \mathbb{R})$ thanks to the singular Adams inequality. A straightforward
calculation shows that
  $$\langle I'(u),  \phi\rangle=M(\|u\|_{E}^{2})\int_{\mathbb{R}^{4}}(\triangle u \triangle \phi+V(x)u \phi )dx-\int_{\mathbb{R}^{4}}\frac{f(x, u)}{|x|^{\eta}}\phi dx, \eqno(3.2)$$
for all\ $u, \phi\in E$,
hence, a critical point of\ (3.2) is a weak solution of\ $(P)$.\\
\textbf{Lemma 3.1} Assume that $(f_{2})$ and $(f_{4})$ hold. Then there exists positive constants $ \delta$ and $r$ such that
\begin{center}
{$I(u)\geq \delta$} for {$\|u\|_{E}=r$}.
\end{center}
$\mathbf{Proof.}$ From $(f_{4})$, there exist $\sigma, \epsilon>0$, such that if $\|u\|_{E}\leq  \epsilon$,
 $$F(x,u)\leq \frac{\mathfrak{M}(1)\lambda_{\eta}-\sigma}{2}|u|^{4},  $$
 for all $x\in \mathbb{R}^{4}$. On the other hand, using $(f_{2})$ for each $q> 4$, we have
\begin{equation*}
\left.
\begin{aligned}[b]
    F(x,u)\leq &\frac{c_{1}}{4}|u|^{4}+c_{2}|u|(e^{\alpha_{0} u^{2}}-1)\\
        \leq & C|u|^{q}(e^{\alpha_{0} |u|^{2}}-1)
\end{aligned}
\right.
\end{equation*}
for $\|u\|_{E}\geq  \epsilon$ and $x\in \mathbb{R}^{4}$. Combining the above estimates, we obtain
\begin{equation*}
\left.
\begin{aligned}[b]
 F(x,u)\leq & \frac{\mathfrak{M}(1)\lambda_{\eta}-\sigma}{2}|u|^{4}+C|u|^{q}(e^{\alpha_{0} |u|^{2}}-1)
\end{aligned}
\right.
\end{equation*}
 for all $(x, u)\in \mathbb{R}^{4}\times \mathbb{R}$.
 On the other hand, \eqref{1.4} gives $\mathfrak{M}(t)\geq \mathfrak{M}(1)t^2, t\in [0, 1]$. Fixed $r>0$ and  $\|u\|_{E}\leq r\leq 1$  such that $\alpha_{0} r^2<32\pi^2(1-\frac{\eta}{4})$, then Lemma 2.2  implies
\begin{equation*}
\left.
\begin{aligned}[b]
            I(u)=&\frac{1}{2}\mathfrak{M}(\|u\|_{E}^{2}) -\int_{\mathbb{R}^{4}}\frac{F(x, u)}{|x|^{\eta}}dx\\
            \geq&\frac{ \mathfrak{M}(1)}{2}\|u\|_{E}^{4}-\frac{\mathfrak{M}(1)\lambda_{\eta}-\sigma}{2}\int_{\mathbb{R}^{4}}\frac{|u|^{4}}{|x|^{\eta}}dx-C\int_{\mathbb{R}^{4}}\frac{|u|^{q}(e^{\alpha_{0} |u|^{2}}-1)}{|x|^{\eta}}dx\\
             \geq&\frac{\mathfrak{M}(1)}{2}\|u\|_{E}^{4}-\frac{\mathfrak{M}(1)\lambda_{\eta}-\sigma}{2}\int_{\mathbb{R}^{4}}\frac{|u|^{4}}{|x|^{\eta}}dx-C\|u\|_{E}^{q}\\
             \geq&\frac{\mathfrak{M}(1)}{2}\|u\|_{E}^{4}-\frac{\mathfrak{M}(1)\lambda_{\eta}-\sigma}{2\lambda_{\eta}}\|u\|_{E}^{4}-C\|u\|_{E}^{q}\\
             =&\frac{\sigma}{2\lambda_{\eta}}\|u\|_{E}^{4}-C\|u\|_{E}^{q}.
\end{aligned}
\right.
\end{equation*}
 Hence, $I$ is bounded form below for $\|u\|_{E}\leq r\leq 1$. Since $\sigma >0$ and $q >4$, we may choose sufficiently small $r>0$ such that
 $$\frac{\sigma}{2\lambda_{\eta}}r^{4}-Cr^{q}\geq \frac{\sigma}{4\lambda_{\eta}}r^{4},$$
 we derive that
 \begin{center}
{$I(u)\geq \frac{\sigma}{4\lambda_{\eta}}r^{4}:=\delta$} for {$\|u\|_{E}=r$.}
\end{center}
This completes the proof. $\hfill\Box$\\
\textbf{Lemma 3.2} Assume $(f_{3})$ is satisfied. Then there exists $e\in B_{r}^{c}(0)$ such that
\begin{center}
{$I(e)< \inf \limits_{\|u\|_{E}=r}I(u)$},
\end{center}
where $r$ are given in Lemma 3.1.\\
$\mathbf{Proof.}$  From $(M_3)$, we have $\mathfrak{M}(t)\leq \mathfrak{M}(1)t^2, t\geq 1$. Let $u\in E\setminus \{0\}$, $u\geq0$ with compact support $\Omega =supp (u)$ and $\|u\|=1$, by $(f_{3})$, for $\mu >4$, there exists $C_1, C_2>0$ such that
for all $(x, s) \in \Omega\times \mathbb{R}^{+}$,
   $$F(x, s)\geq C_1 s^{\mu}-C_2.$$
Then
\begin{equation*}
\left.
\begin{aligned}[b]
 I(tu)\leq& \frac{\mathfrak{M}(1)t^{4}}{2}\|u\|_{E}^{4}-C_1t^{\mu}\int_{\Omega}\frac{|u|^{\mu}}{|x|^{\eta}}dx+C_2|\Omega|,
\end{aligned}
\right.
\end{equation*}
which implies that $I(tu)\rightarrow -\infty$ as $t\rightarrow \infty$. Setting $e = tu$ with $t$ sufficiently large, we finish the proof of the lemma.  $\hfill\Box$

From Lemma 3.1, Lemma 3.2,   we get  a $(PS)_{c}$ sequence $\{u_{n}\}\subset E$, i.e.
\be I(u_{n})\rightarrow c>0  ~~\operatorname{and}~~ I'(u_{n})\rightarrow 0 ~~~~~\operatorname{as}~~ n \rightarrow\infty,   \lb{(3.3)}\ee
where
\be c=\inf \limits _{\gamma\in \Gamma} \max \limits _{t\in [0, 1]} I(\gamma(t))\label{999}\ee
and
$$\Gamma=:\{\gamma\in C([0, 1]:E):\gamma(0)=0, \gamma(1)=e\}.$$
\textbf{Lemma 3.3} Suppose $(f_{6})$ is satisfied, then the level $c\in\bigg(0, \frac{1}{2}\mathfrak{M}(\frac{32\pi^{2}}{\alpha_{0}}\left(1-\frac{\eta}{4}\right))\bigg)$.\\
\textbf{Proof.} Firstly, we claim the best constant $S_{p}$ can be obtained. In fact, since
$$S_{p}=\inf_{u\in E\backslash\{0\}}\frac{\| u \|_{E}}{(\int_{\mathbb{R}^{4}}\frac { |u|^{p}}{|x|^{\eta}}dx)^{\frac{1}{p}}},$$
we can choose $u_{n}$ such that
   $$\int_{\mathbb{R}^{4}}\frac { |u_{n}|^{p}}{|x|^{\eta}}dx=1  ~~\operatorname{and}~~ \| u_{n} \|_{E}\rightarrow S_{p}~~\operatorname{as}~~ n \rightarrow\infty,$$
so $u_{n}$ is bounded in $E$. From Lemma 2.1, there exists $u_0\in E$ such that up to a subsequence $u_{n}\rightharpoonup u_0$ in $E$,
$u_{n}\rightarrow u_0$ in $L^{p}(\mathbb{R}^{4}, |x|^{-\eta}dx)$ and $u_{n}(x)\rightarrow u_0(x)$ almost everywhere in $\mathbb{R}^{4}$.
This implies
  $$\int_{\mathbb{R}^{4}}\frac { |u_0|^{p}}{|x|^{\eta}}dx=\lim_{n\rightarrow \infty}\int_{\mathbb{R}^{4}}\frac { |u_{n}|^{p}}{|x|^{\eta}}dx=1.$$
We also have $\| u_0 \|_{E}\leq\lim_{n\rightarrow \infty}\| u_{n} \|_{E}=S_{p}$, thus $\| u_0 \|_{E}=S_{p}$. From the definition of $c$, let $\gamma:[0, 1]\rightarrow E, \gamma(t)=tt_{0}u$, where\ $t_{0}$ is a real number which satisfies\ $I(t_{0}u_0)<0$, we have\ $\gamma\in \Gamma$, and
therefore
   $$c\leq\max_{t\in[0, 1]}I(\gamma(t))\leq\max_{t\geq0}I(tu_0)=\max_{t\geq 0}\bigg(\frac{\mathfrak{M}(t^{2}S_{p}^{2})}{2}-\int_{\mathbb{R}^{4}}\frac{F(x, tu_0)}{|x|^{\eta}}dx  \bigg),$$
by $(f'_{6})$, we have
   $$c\leq \max_{t\geq 0}I(tu)=\max_{t\geq 0}\bigg(\frac{\mathfrak{M}(t^{2}S_{p}^{2})}{2}-\frac{t^{p}}{p}C_{p}\bigg)
   <\frac{1}{2}\mathfrak{M}(\frac{32\pi^{2}}{\alpha_{0}}\left(1-\frac{\eta}{4}\right)).$$
The proof of the lemma is completed. $\hfill\Box$

Consider the Nehari manifold associated to the functional $I$, that is,
  $$\mathcal{N}:=\{u\in E\backslash\{0\}:I'(u)u=0\}$$
and $c^*=\inf_{u\in \mathcal{N}}I(u).$\\
\textbf{Lemma 3.4}  Suppose $M$ satisfies $(M_{3})$, $f$ satisfies $(f_{5})$. Then $c\leq c^*.$\\
\textbf{Proof.}  Let\ $u\in \mathcal{N}$, we define\ $h:(0, +\infty)\rightarrow \mathbb{R}$ by\ $h(t)=I(tu)$.
 We have that $h$ is differentiable and
                         $$h'(t)=I'(tu)u=M(t^{2}\|u\|^{2})t\|u\|^{2}-\int_{\mathbb{R}^{4}}\frac{f(x, tu)u}{|x|^{\eta}}dx,~~~~\forall t\geq0.$$
From\ $I'(u)u=0$, we get
 $$h'(t)=I'(tu)u-t^{3}I'(u)u,$$
so
\begin{equation*}
 \left.
   \begin{aligned}[b]
    h'(t)=& t^{3}\|u\|_E^{4}\bigg[\frac{M(t^{2}\|u\|^{2})}{t^{2}\|u\|_E^{2}}-\frac{M(\|u\|^{2})}{\|u\|_E^{2}}\bigg]\\
     &+  t^{3}\int_{\mathbb{R}^{4}}\bigg(\frac{f(x, u)}{u^{3}}-\frac{f(x, tu)}{(tu)^{3}}\bigg)u^{4}dx,
\end{aligned}
\right.
\end{equation*}
By $(M_3)$, $(f_{5})$, we conclude that $h'(t)>0$ for $0< t<1$ and $h'(t)<0$ for $t>1$. Thus, $h(1)=\max_{t\geq0}h(t)$, which means
                                             $$I(u)=\max_{t\geq0}I(tu).$$
From the above argument, we see that $h'(t)<0$ is strongly decreasing in $t\in (1,+\infty)$, so $h(t)\to -\infty$ as $t\to +\infty$.
Now, define\ $\gamma:[0, 1]\rightarrow E, \gamma(t)=tt_{0}u$, where\ $t_{0}$ is a real number which satisfies\ $I(t_{0}u)<0$, we have\ $\gamma\in \Gamma$, and
therefore
                      $$c\leq\max_{t\in[0, 1]}I(\gamma(t))\leq\max_{t\geq0}I(tu)=I(u).$$
Since\ $u\in  \mathcal{N}$ is arbitrary, we have $c\leq c^*$.  $\hfill\Box$
\section{The ground state solution}
In this section, we consider the ground state solution. We first prove
the following convergence results.\\
\textbf{Lemma 4.1} Suppose $(V_1)$, $(f_{1})-(f_{5})$  are satisfied, let\ $\{u_{n}\}$ is an arbitrary $(PS)_{c}$ sequence, then there exists a subsequence of\ $\{u_{n}\}$(still denoted by $\{u_{n}\}$) and $u\in E$ such that
\begin{equation*}
 \left\{  \begin{array}{l}
          \frac{f(x, u_{n})}{|x|^{\eta}}\rightarrow \frac{f(x, u)}{|x|^{\eta}}~~~~strongly~~in~~L^{1}_{loc}(\mathbb{R}^{4}), \\
           \frac{F(x, u_{n})}{|x|^{\eta}}\rightarrow \frac{F(x, u)}{|x|^{\eta}}~~~~strongly~~in~~L^{1}(\mathbb{R}^{4}), \\
         \end{array}
         \right.
 \end{equation*}
\textbf{Proof.}  Let\ $\{u_{n}\}\subset E$ be an arbitrary $(PS)_{c}$ sequence of $I$, {\it i.e.}
\be I(u_{n})\rightarrow c>0  ~~\operatorname{and}~~ I'(u_{n})\rightarrow 0 \quad \operatorname{as}~~ n \rightarrow\infty.  \lb{4.1}\ee
We shall prove that the sequence $\{u_{n}\}$ is bounded in $E$. Indeed, since $\mu>4$,
then
 \begin{equation*}
\left.
\begin{aligned}[b]
            c+o_{n}(1)\|u_n\|_E\geq &I(u_{n})-\frac{1}{\mu}\langle I'(u_n), u_n\rangle\\
            \geq &\frac{1}{2}\mathfrak{M}(\|u_n\|_E^2)-\frac{1}{\mu}M(\|u_n\|_E^2)\|u_n\|_E^2-\frac{1}{\mu}\int_{\mathbb{R}^{4}}\frac{\mu F(x, u_n)-f(x, u_n)u_n}{|x|^\eta}dx\\
            \geq &(\frac{1}{4}-\frac{1}{\mu})M(\|u_n\|_E^2)\|u_n\|_E^2-\frac{1}{\mu}\int_{\mathbb{R}^{4}}\frac{\mu F(x, u_n)-f(x, u_n)u_n}{|x|^\eta}dx\\
            \geq &(\frac{1}{4}-\frac{1}{\mu})M_0\|u_n\|_E^2,
\end{aligned}
\right.
\end{equation*}
which implies that\ $\{u_{n}\}$ is bounded in $E$.
It then follows from \eqref{4.1} that
$$\frac{f(x, u_{n})u_{n}}{|x|^{\eta}}dx\leq C, ~~~~~\frac{F(x, u_{n})}{|x|^{\eta}}dx\leq C.$$
By Lemma 2.1 of \cite{19}, we get
\be \frac{f(x, u_{n})}{|x|^{\eta}}\rightarrow \frac{f(x, u)}{|x|^{\eta}}~~~~strongly~~in~~L^{1}_{loc}(\mathbb{R}^{N}).\lb{4.2}\ee
By $(f_{2})$ and $(f_{3})$, there exists $C>0$ such that
  $$F(x, u_{n})\leq C_1|u_{n}|^{4}+C_2f(x, u_{n}).$$
From Lemma 2.2 and generalized Lebesgue's dominated convergence theorem,  arguing as Lemma 4.7 in \cite{L. Chen}, we can derive that
          \be \frac{F(x, u_{n})}{|x|^{\eta}}\rightarrow \frac{F(x, u)}{|x|^{\eta}}~~~~strongly~~in~~L^{1}(\mathbb{R}^{4}). \lb{4.3}\ee
This completes the proof of the lemma.  $\hfill\Box$\\
\textbf{Lemma 4.2} Let $(M_1)-(M_3)$ and $(f_1)-(f_6)$ hold. Then the functional $I$ satisfies the $(PS)_c$ condition.\\
\textbf{Proof.} By the process in proof of Lemma 4.1, we have that the $(PS)_{c}$ sequence $\{u_{n}\}$ is bounded in\ $E$. We claim that $I(u)\geq 0 $. Indeed, suppose by contradiction that $I(u) < 0$. Then $u\neq 0$, set $r(t):=I(tu), t\geq 0$, we have $r(0)=0$ and $r(1)<0$. As the proof of Lemma 3.1, for $t>0$ small enough, it holds $r(t)>0$. So there exists $t_0\in (0, 1)$  such that
    $$r(t_0)=\max_{t\in [0, 1]}r(t),~~~~r'(t_0)=\langle I'(t_0u), u\rangle=0,$$
By Remark 1.1 and Lemma 2.3, we have
\begin{equation*}
\left.
\begin{aligned}[b]
         c\leq I(t_0u)=&I(t_0u)-\frac{1}{4}  \langle I'(t_0 u), u   \rangle\\
            =&\frac{1}{2}\mathfrak{M}(\|t_{0}u\|_E^{2})-\frac{1}{4}M(\|t_{0}u\|_E^{2})\|t_{0}u\|_E^{2}\\
             &+\frac{1}{4}\int_{\mathbb{R}^{4}}\frac{f (x, t_0u)t_0u-4F (x, t_0u)}{|x|^{\eta}}dx\\
             <&\frac{1}{2}\mathfrak{M}(\|u\|_E^{2})-\frac{1}{4}M(\|u\|_E^{2})\|u\|_E^{2}\\
              &+\frac{1}{4}\int_{\mathbb{R}^{4}}\frac{f (x, u)u-4F (x, u)}{|x|^{\eta}}dx
\end{aligned}
\right.
\end{equation*}
FUrthmore, by the weak lower semicontinuity of the norm and Fatou's Lemma, we have
\begin{equation*}
\left.
\begin{aligned}[b]
         c<&\liminf_{n\rightarrow \infty}(\frac{1}{2}\mathfrak{M}(\|u_n\|_E^{2})-\frac{1}{4}M(\|u_n\|_E^{2})\|u\|_E^{2})\\
            &+\frac{1}{4}\liminf_{n\rightarrow \infty}\int_{\mathbb{R}^{4}}\frac{f (x, u_n)u_n-4F (x, u_n)}{|x|^{\eta}}dx\\
             \leq& \liminf_{n\rightarrow \infty}(I'(u_n)-\frac{1}{4}\langle I'(u_n), u_n   \rangle)=c.\\
\end{aligned}
\right.
\end{equation*}
which is not impossible. Thus the claim is true. From
the lower semi-continuity of the norm in $E$, we have $\|u\|_E\leq\lim_{n\rightarrow\infty}\|u_n\|_E$.  Suppose,
by contradiction, that $\|u\|_E<\lim_{n\rightarrow\infty}\|u_n\|_E:=\xi$. Set $v_n:=\frac{u_n}{\|u_n\|_E}$ and $v:=\frac{u}{\xi}$, then $v_n\rightharpoonup v $ weakly in $E$ and $\|v\|_E<1$. From $I(u)\geq 0$ and Lemma 4.1, we have
\begin{equation*}
\left.
\begin{aligned}[b]
         \mathfrak{M}(\xi^2)=&\lim_{n\rightarrow\infty}\mathfrak{M}(\|u_n\|_E^2)=\lim_{n\rightarrow\infty}(2(I(u_n)+\int_{\mathbb{R}^{4}}\frac{F(x, u_n)}{|x|^{\eta}}dx)\\
         =&2c+2\int_{\mathbb{R}^{4}}\frac{F(x, u)}{|x|^{\eta}}dx=2c+\mathfrak{M}(\|u\|_E^2)-2I(u)\\
         <&\mathfrak{M}((1-\frac{\eta}{4})\frac{32\pi ^2}{\alpha_{0}})+\mathfrak{M}(\|u\|_E^2)\\
         \leq&\mathfrak{M}\bigg((1-\frac{\eta}{4})\frac{32\pi ^2}{\alpha_{0}})+\|u\|_E^2\bigg)
\end{aligned}
\right.
\end{equation*}
Here, we have used the condition $(M_2)$ in the last inequality. Since $\mathfrak{M}$ is increasing, it holds $\xi^2< (1-\frac{\eta}{4})\frac{32\pi ^2}{\alpha_{0}}+\|u\|_E^2$. Notice that
   $$\xi^2=\frac{\xi^2-\|u\|_E^2}{1-\|v\|_E^2}.$$
Thus
  $$\xi^2=\frac{(1-\frac{\eta}{4})\frac{32\pi ^2}{\alpha_{0}}}{1-\|v\|_E^2}.$$
Choosing $q>1$ sufficiently close to 1 and $\beta_0>0$ such that for large $n$,
      $$q\alpha_0\|u_n\|_E^{2}\leq \beta_0<\frac{(1-\frac{\eta}{4})32\pi ^2}{1-\|v\|_E^2}.$$
From concentration compactness principle with singular Adams inequality, we have
\begin{equation}
\left.
\begin{aligned}[b]
\int_{\mathbb{R}^{4}}\frac{e^{q\alpha_{0} u_{n}^{2}-1}}{|x|^{\eta}}dx\leq \int_{\mathbb{R}^{4}}\frac{e^{\beta_{0} v_{n}^{2}-1}}{|x|^{\eta}}dx\leq C.
\end{aligned}
\right. \lb{(6.10)}
\end{equation}
From $(f_{2})$ and H\"{o}lder inequality, we have
\begin{equation}
\left.
\begin{aligned}[b]
&\left|\int_{\mathbb{R}^{4}}\frac{f(x, u_{n})(u_{n}-u)}{|x|^{\eta}} dx\right|\\
&\leq c_{1}\bigg(\int_{\mathbb{R}^{4}}\frac{|u_{n}|^{4}}{|x|^{\eta}}dx\bigg)^{\frac{3}{4}}\bigg(\int_{\mathbb{R}^{4}}\frac{|u_{n}-u|^{4}}{|x|^{\eta}}dx\bigg)^{\frac{1}{4}}\\
&+c_{2}\bigg(\int_{\mathbb{R}^{4}}\frac{|u_{n}-u|^{q'}}{|x|^{\eta}}dx\bigg)^{\frac{1}{q'}} \bigg(\int_{\mathbb{R}^{4}}\frac{e^{q\alpha_{0} u_{n}^{2}-1}}{|x|^{\eta}}dx \bigg)^{\frac{1}{q}},
\end{aligned}
\right.
\lb{(6.11)}
\end{equation}
where $\frac{1}{q'}+\frac{1}{q}=1$. In view Lemma 2.1, combining  \eqref{(6.10)}  with   \eqref{(6.11)}, we obtain
\be \int_{\mathbb{R}^{N}}\frac{f(x, u_{n})(u_{n}-u)}{|x|^{\eta}}dx\rightarrow 0.  \lb{(6.12)}\ee
Since $I'(u_n)(u_n -u) \rightarrow 0$, we have
\be M(\|u_n\|_E^2)\int_{\mathbb{R}^{4}}(\Delta u_{n} \Delta (u_{n}-u) +V(x)u_{n}(u_{n}-u) )dx\rightarrow 0, \lb{(6.13)}\ee
On the other hand, by $u_{n}\rightharpoonup u$ in $E$, we have
\be M(\|u_n\|_E^2)\int_{\mathbb{R}^{4}}\big(\Delta u \Delta(u_{n}-u) +V(x)u(u_{n}-u)\big)dx \rightarrow 0. \lb{(6.14)}\ee
\eqref{(6.13)} minus \eqref{(6.14)},
 we can derive
\begin{equation}
\left.
\begin{aligned}[b]
 &\lim_{n\rightarrow \infty}M(\|u_n\|_E^2)\|u_{n}-u\|_{E}^{2}=0,
  \end{aligned}
\right.
\lb{(6.12)}
\end{equation}
which is in contradiction with the fact $\|u\|_E<\lim_{n\rightarrow\infty}\|u_n\|_E:=\xi$. Thus, we have $\|u\|_E=\xi=\lim_{n\rightarrow\infty}\|u_n\|_E$. Since $\{u_n\}$ is bounded in $E$, we can apply Brezis-Lieb lemma to obtain $u_n\rightarrow u$ strongly in $E$. $\hfill\Box$\\
\noindent\textbf{The proof of Theorem 1.1.} Since $I\in C^1(E, \mathbb{R})$, by Lemma 4.2, we have  $I'(u) = 0$ and $I(u) = c$. Therefore, by the definition of $c^*$ and $c\leq c^*$, we know $u$ is a ground state solution.

Next, we will show that\ $u$ is nonzero. If\ $u\equiv0$, since $F(x, 0)=0$ for all $x\in \mathbb{R}^{4}$, from Lemma 3.4, we have
                          \be\lim_{n\rightarrow\infty}\frac{1}{2}\mathfrak{M}(\|u_{n}\|_{E}^{2})c<\mathfrak{M}\bigg((1-\frac{\eta}{4})\frac{32\pi^2}{\alpha_{0}}\bigg), \lb{(4.1)}\ee
Thus, there exist some $\epsilon_{0}>0$ and $n_{*}>0$ such that $\|u_{n}\|_{E}^{2}\leq (1-\frac{\eta}{4})\frac{32\pi^2}{\alpha_{0}}-\epsilon_{0}$ for all $n>n_{*}$. Choose $q>1$ sufficiently close to 1 such that $q\alpha_{0}\|u_{n}\|_{E}^{2}\leq (1-\eta /4)32\pi^2-\epsilon_{0}\alpha_{0}$ for all $n>n_{*}$.
By $(f_{2})$, there holds
    $$|f(x, u_{n})u_n|\leq c_{1}|u_{n}|^{\
    4}+c_{2}|u_{n}|(e^{\alpha_{0} u_{n}^{2}}-1).$$
Thus by using singular Adams inequality, we have
\begin{equation*}
\left.
\begin{aligned}[b]
         &\int_{\mathbb{R}^{4}}\frac{|f(x, u_{n})u_{n}|}{|x|^{\eta}}dx\\
         \leq&c_{1}\int_{\mathbb{R}^{4}}\frac{|u_{n}|^{4}}{|x|^{\eta}}dx+c_{2}\int_{\mathbb{R}^{4}}\frac{|u_{n}|(e^{\alpha_{0} |u_{n}|^{2}}-1)}{|x|^{\eta}}dx\\
         \leq&c_{1}\int_{\mathbb{R}^{4}}\frac{|u_{n}|^{4}}{|x|^{\eta}}dx+c_{2}\bigg(\int_{\mathbb{R}^{4}}\frac{e^{q\alpha_{0} |u_{n}|^{2}}-1}{|x|^{\eta}}dx\bigg)^{\frac{1}{q}}
         \left(\int_{\mathbb{R}^{4}}\frac{|u_{n}|^{q'}}{|x|^{\eta}}dx\right)^{\frac{1}{q'}}\\
         \leq&c_{1}\int_{\mathbb{R}^{4}}\frac{|u_{n}|^{4}}{|x|^{\eta}}dx+C
         \left(\int_{\mathbb{R}^{4}}\frac{|u_{n}|^{q'}}{|x|^{\eta}}dx\right)^{\frac{1}{q'}}\rightarrow 0,
\end{aligned}
\right.
\end{equation*}
here we have used Lemma 2.1 in the last estimate. From $I^{'}(u_{n})u_{n}\rightarrow 0$, we have
\be\lim_{n\rightarrow\infty}M(\|u_{n}\|_{E}^{2})\|u_{n}\|_{E}^{2}=0,  \lb{(4.2)}\ee
From the condition $(M_1)$, we can get $\|u_n\|\rightarrow 0$. Then $I(u_n)\rightarrow 0$, which contradics the fact that  $I(u_n)\rightarrow c>0$, so $u$ is nonzero. From $I (u) = c > 0$, we know $u$ is positive.
This completes the proof of Theorem 1.1. ~~~~~~~~$\hfill\Box$
\section{The ground state  solution without the A-R condition}
In this section, we instead the condition  $(f_3)$,  the nonlinear term satisfies the exponential growth but without satisfying the Ambrosetti-Rabinowitz condition, we assume that

$(f'_{3})$ $\lim_{|t|\rightarrow +\infty}\frac{F(x, t)}{|t|^{4}}=\infty$ uniformly on $x \in \mathbb{R}^{4}$, where $F(x,t)=\int_0^tf(x,s)ds$.\\
 We will use a Cerami's Mountain Pass Theorem which was introduced in \cite{16a, 17}. The detail is the following:\\
\textbf{Definition A.} Let $(E, \|\cdot\|_{E})$ be a real Banach space with its dual space $(E^{*}, \|\cdot\|_{E^{*}})$. Suppose $I\in C^{1}(E, \mathbb{R})$. For  $c\in \mathbb R$, we say that  $\{u_{n}\}\subset E$   a $(C)_{c}$ sequence of the functional $I$, if
$$I(u_{n})\rightarrow c  ~~\operatorname{and}~~ (1+\|u_{n}\|_{E})\|I'(u_{n})\|_{E^{*}}\rightarrow 0 ~~~~~\operatorname{as}~~ n \rightarrow\infty.   $$
\textbf{Proposition A.}  Let $(E, \|\cdot\|_{E})$ be a real Banach space, $I\in C^{1}(E, \mathbb{R})$,  $I(0)=0$ and satisfies:

$(i)$ there exists positive constants $ \delta$ and $r$ such that
\begin{center}
{$I(u)\geq \delta$} for {$\|u\|_{E}=r$}
\end{center}
and

 $(ii)$ there exists $e\in E$ with $\|e\|_E>r$ such that
\begin{center}
{$I(e)\leq 0$}.
\end{center}

Define $c$  by
$$c=\inf \limits _{\gamma\in \Gamma} \max \limits _{t\in [0, 1]} I(\gamma(t)),$$
where
$$\Gamma=:\{\gamma\in C([0, 1]:E):\gamma(0)=0, \gamma(1)=e\}.$$
Then $I$ possesses a $(C)_{c}$ sequence.

Firstly, we check the geometry of the functional $I$ under the weak condition. Secondly, the key to establish the results in previous sections is prove that the Cerami sequence is bounded. Once we will have proved this, the remaining parts are similar.  \\
\textbf{Lemma 5.1.}  Assume that $(V_{1})$, $(f_{2})$-$(f_{4})$ hold. Then

 (i) there exists positive constants $ \delta$ and $r$ such that
\begin{center}
{$I(u)\geq \delta$} for {$\|u\|_{E}=r$}.
\end{center}

(ii) there exists $e\in E$ with $\|e\|_E>r$ such that
\begin{center}
{$I(e)< \inf \limits_{\|u\|_{E}=r}I(u)$},
\end{center}
$\mathbf{Proof.}$  The proof of (i) is similar as Lemma 3.1. From $(M_3)$, we have $\mathfrak{M}(t)\leq \mathfrak{M}(1)t^2, t\geq 1$. Let $u\in E\setminus \{0\}$, $u\geq0$ with compact support $\Omega =supp (u)$, by $(f_{2})$, for all $L$, there exists $d$ such that
for all $(x, s) \in \Omega\times \mathbb{R}^{+}$,
   $$F(x, s)\geq L s^{4}-d.$$
Then
\begin{equation*}
\left.
\begin{aligned}[b]
 I(tu)\leq& \frac{\mathfrak{M}(1) t^{4}}{2}\|u\|_{E}^{4}-Lt^{4}\int_{\Omega}\frac{|u|^{4}}{|x|^{\eta}}dx+O(1)\\
      \leq& t^{4}\bigg(\frac{\mathfrak{M}(1)\|u\|_{E}^{4}}{2}-L\int_{\Omega}\frac{|u|^{4}}{|x|^{\eta}}dx\bigg)+O(1).
\end{aligned}
\right.
\end{equation*}
Now  choosing $L> \frac{\mathfrak{M}(1)\|u\|_{E}^{4}}{N\int_{\Omega}\frac{|u|^{4}}{|x|^{\eta}}dx}$, it implies that $I(tu)\rightarrow -\infty$ as $t\rightarrow \infty$.   Setting $e = tu$ with $t$ sufficiently large, the proof of
(ii) is completed.  $\hfill\Box$

From Lemma 3.1, Lemma 5.1 and Proposition A,   we get  a $(C)_{c}$ sequence $\{u_{n}\}\subset E$, i.e.
\be I(u_{n})\rightarrow c>0  ~~\operatorname{and}~~ (1+\|u_{n}\|_{E})\|I'(u_{n})\|_{E^{*}}\rightarrow 0 ~~~~~\operatorname{as}~~ n \rightarrow\infty,   \lb{(5.1)}\ee
\textbf{Lemma 5.2.}  Let\ $\{u_{n}\}\subset E$ be an arbitrary Cerami sequence of $I$, Then $\{u_n\}$ is bounded up to a subsequence.\\
\textbf{Proof.} Let\ $\{u_{n}\}\subset E$ be an arbitrary Cerami sequence of $I$, {\it i.e.}
 \be\frac{\mathfrak{M}(\|u_n\|_{E}^{2})}{2} -\int_{\mathbb{R}^{4}}\frac{F(x, u_n)}{|x|^{\eta}}dx\rightarrow c   ~~~\operatorname{as}~~ n \rightarrow\infty, \lb{(aaa)}\ee
and
\be (1+\|u_{n}\|_{E})|\langle I'(u_{n}), \varphi \rangle|\leq \tau_n \|\varphi\|_{E}~~~~~\operatorname{for ~~all}~~ \varphi\in E,  \lb{(bbb)} \ee
where $\tau_n\rightarrow 0$ as $n\rightarrow \infty$. We shall prove that the sequence $\{u_{n}\}$ is bounded in $E$. Indeed, suppose by contradiction that
$$ \|u_{n}\|_E\rightarrow +\infty$$
and set
$$v_{n}=\frac{u_{n}}{\|u_{n}\|_E},$$
then $\|v_{n}\|=1$. From Lemma 2.1, we can assume that for any $q\geq 4$, there exists $v\in E$ such that up to a subsequence
\begin{equation*}
 \left\{  \begin{array}{l}
   v_{n}^{+}\rightharpoonup v^{+} ~~\operatorname{in}~~ E,\\
     v_{n}^{+}\rightarrow v^{+} ~~\operatorname{in}~~ L^{q}(\mathbb{R}^{4}),\\
    v_{n}^{+}\rightarrow v^{+} ~~\operatorname{a.e.~in}~~\mathbb{R}^{4}.
\end{array}
  \right.
\end{equation*}
We will show that $v^{+}=0$ a.e. in $\mathbb{R}^{4}$. In fact, if $ \Lambda^{+}=\{x\in \mathbb{R}^{4}: v^{+}(x)>0\}$ has a positive measure, then in $\Lambda^{+}$, we have
    $$\lim_{n\rightarrow\infty}u_{n}^{+}=\lim_{n\rightarrow\infty}v_{n}^{+}\|u_{n}\|=+\infty.$$
From $(f'_{3})$ we have
   $$\lim_{n\rightarrow \infty}\frac{F(x, u_{n}^{+}(x))}{|x|^\eta|u_{n}^{+}(x)|^{4}}=+\infty~~\operatorname{a.e.~in}~~\Lambda^{+},$$
and
  $$\lim_{n\rightarrow \infty}\frac{F(x, u_{n}^{+}(x))}{|x|^\eta|u_{n}^{+}(x)|^{4}}|v_{n}^{+}(x)|^{4}=+\infty~~\operatorname{a.e.~in}~~\Lambda^{+}.$$
Thus
    $$\int\limits_{\mathbb{R}^{4}}\liminf_{n\rightarrow \infty}\frac{F(x, u_{n}^{+}(x))}{|x|^\eta|u_{n}^{+}(x)|^{4}}|v_{n}^{+}(x)|^{4}dx=+\infty.$$
Since $\{u_{n}\}\subset E$ be an arbitrary Cerami sequence of $I$, we have
     $$\mathfrak{M}(\|u_{n}\|_{E}^{2})=2c+2\int\limits_{\mathbb{R}^{4}}\frac{F(x, u_{n}^{+}(x))}{|x|^\eta}dx+o_{n}(1)$$
 Since $\mathfrak{M}$ is increasing, it holds
   $$\int\limits_{\mathbb{R}^{4}}\frac{F(x, u_{n}^{+}(x))}{|x|^\eta}dx\rightarrow +\infty.$$
From $(M_3)$, we have $\mathfrak{M}(t)\leq \mathfrak{M}(1)t^2, t\geq 1$. Thus
\begin{equation*}
\left.
\begin{aligned}[b]
  &\liminf_{n\rightarrow \infty}\int\limits_{\mathbb{R}^{4}}\frac{F(x, u_{n}^{+}(x))}{|x|^\eta|u_{n}^{+}(x)|^{4}}|v_{n}^{+}(x)|^{4}dx\\
  & =\liminf_{n\rightarrow \infty}\int\limits_{\mathbb{R}^{4}}\frac{F(x, u_{n}^{+}(x))}{|x|^\eta\|u_{n}\|^{4}}dx\\
   & \leq\liminf_{n\rightarrow \infty}\int\limits_{\mathbb{R}^{4}}\frac{\mathfrak{M}(1)F(x, u_{n}^{+}(x))}{|x|^\eta\mathfrak{M}(\|u_{n}\|_{E}^{2})}dx\\
  & =\liminf_{n\rightarrow \infty} \frac{\int\limits_{\mathbb{R}^{4}}\frac{F(x, u_{n}^{+}(x))}{|x|^\eta}dx}{2c+2\int\limits_{\mathbb{R}^{4}}\frac{F(x, u_{n}^{+}(x))}{a(x)}dx+o_{n}(1)}\\
  &=\frac{1}{2}.
\end{aligned}
\right.
\end{equation*}
This is  a contradiction. Hence $v\leq 0$ a.e. and  $v_n^{+}\rightharpoonup 0$ in $E$.

Let $t_{n}\in [0, 1]$ be such that
 $$I(t_{n}u_{n})=\max_{t\in[0, 1]}I(tu_{n}).$$
For any given $A\in \big(0, \left((1-\frac{\eta}{4})\frac{32\pi^2}{\alpha_{0}}\right)^\frac{1}{2}\big)$, for the sake of simplicity, let $$\epsilon=\frac{(1-\frac{\eta}{4})32\pi ^2}{A^{2}}-\alpha_{0}>0.$$
In the following argument we will take $A\to \left((1-\frac{\eta}{4})\frac{32\pi^2}{\alpha_{0}}\right)^\frac{1}{2}$ and so we have $\epsilon\to 0$.\\
By condition $(f_{2})$,  there
exists $C > 0$ such that
 \be F(x, t)\leq C|t|^{4}+\epsilon R(\alpha_{0}+\epsilon, |t|),~~~~\forall (x, t)\in  \mathbb{R}^{4}\times \mathbb{R}^{+},\lb{aa}\ee
where $R(\alpha, s)=e^{\alpha s^{2}}-1$. In fact, from condition $(f_{2})$,  there holds
$$F(x,t)\le \frac CN|t|^{4}+|t|R(\alpha_0,|t|).$$
By using Young inequality, for $\frac 1p+\frac 1q=1$, $p,q>1$, there holds
$$ab \le \epsilon \frac{a^p}{p}+\epsilon^{-q/p} \frac{b^q}{q}.$$
So we have
$$F(x,t)\le \frac CN|t|^{4}+\frac{\epsilon R(\alpha_0,|t|)^p }{p}+\epsilon^{-q/p} \frac{|t|^q}{q}.$$
Now we  take $p=\frac{\alpha_0+\epsilon}{\alpha_0}$ and $q=\frac {\alpha_0+\epsilon}{\epsilon}>4$.    One can see that near infinity $|t|^q$ can be estimated from above by $R(\alpha_{0}+\epsilon, |t|)$, and near the origin  $|t|^q$ can be estimated from above by $|t|^{4}$, thus we obtain \eqref{aa}. We also have $\frac{A}{\|u_{n}\|}\in (0,1]$ with sufficient large $n$, so by using \eqref{aa}, we have
 \begin{equation*}
\left.
\begin{aligned}[b]
            I(t_{n}u_{n})\geq &I(\frac{A}{\|u_{n}\|}u_{n})=I(Av_{n}) =\frac{\mathfrak{M}(A^{2})}{2}-\int_{\mathbb{R}^{4}}\frac{F(x, Av_{n})}{|x|^{\eta}}dx\\
               =&\frac{\mathfrak{M}(A^{N})}{N}-\int_{\mathbb{R}^{4}}\frac{F(x, Av_{n}^{+})}{|x|^{\eta}}dx\\
              \geq &\frac{\mathfrak{M}(A^{2})}{N}-CA^{4}\int_{\mathbb{R}^{4}}\frac{|v_{n}^{+}|^{4}}{|x|^{\eta}}dx
              -\epsilon\int_{\mathbb{R}^{4}}\frac{R(\alpha_{0}+\epsilon, Av_{n}^{+})}{|x|^{\eta}}dx\\
              \geq &\frac{\mathfrak{M}(A^{2})}{2}-CA^{4}\int_{\mathbb{R}^{4}}\frac{|v_{n}^{+}|^{4}}{|x|^{\eta}}dx
              -\epsilon\int_{\mathbb{R}^{4}}\frac{R((\alpha_{0}+\epsilon) A^{2}, v_{n}^{+})}{|x|^{\eta}}dx\\
              \geq &\frac{\mathfrak{M}(A^{2})}{2}-CA^{4}\int_{\mathbb{R}^{4}}\frac{|v_{n}^{+}|^{4}}{|x|^{\eta}}dx
              -\epsilon\int_{\mathbb{R}^{4}}\frac{R((1-\frac{\eta}{4})32\pi^2, v_{n}^{+})}{|x|^{\eta}}dx.\\
\end{aligned}
\right.
\end{equation*}
 Since $v_{n}^{+}\rightharpoonup 0$ in $E$ and the embedding $E\hookrightarrow L^{q}(\mathbb{R}^{4}, |x|^{-\eta}dx)(q\geq 4)$ is compact,
by using the H$\ddot{o}$lder inequality, we have  $\int_{\mathbb{R}^{4}}\frac{|v_{n}^{+}|^{2}}{|x|^{\eta}}dx \rightarrow 0$. By singular Trudinger-Moser inequality,   $\int_{\mathbb{R}^{4}}\frac{R((1-\frac{\eta}{4})32\pi^2, v_{n}^{+})}{|x|^{\eta}}dx$ is bounded. When $  A\rightarrow\big((1-\frac{\eta}{4})\frac{32\pi^2}{\alpha_{0}}\big)^\frac{1}{2}$, we can show
 \be\liminf_{n\rightarrow \infty}I(t_{n}u_{n})\geq \frac{1}{2}\mathfrak{M}\bigg(((1-\frac{\eta}{4})\frac{32\pi^2}{\alpha_{0}})\bigg) >c.  \lb{(3.4)}\ee

Since $I(0)=0$ and $I(u_{n})\rightarrow c$, we can assume $t_{n}\in (0, 1)$, and so $I'(t_{n}u_{n})t_{n}u_{n}=0$, it follows from $(f_{5})$,
\begin{equation*}
\left.
\begin{aligned}[b]
            4I(t_{n}u_{n})=&4I(t_{n}u_{n})- I'(t_{n}u_{n})t_{n}u_{n}\\
             =&2 \mathfrak{M}(\|t_{n}u_{n}\|^{2})-4\int_{\mathbb{R}^{4}}\frac{F(x, t_{n}u_{n})}{|x|^{\eta}}dx\\
             &-M(\|t_{n}u_{n}\|^{2})\|t_{n}u_{n}\|^{2}+
             \int_{\mathbb{R}^{4}}\frac{f(x, t_{n}u_{n})t_{n}u_{n}}{|x|^{\eta}}dx\\
             =&2\mathfrak{M}(\|t_{n}u_{n}\|^{2})-M(\|t_{n}u_{n}\|^{2})\|t_{n}u_{n}\|^{2}+\int_{\mathbb{R}^{4}}\frac{H(x, t_{n}u_{n})}{|x|^{\eta}}dx\\
             \leq &2\mathfrak{M}(\|u_{n}\|^{2})-M(\|u_{n}\|^{2})\|t_{n}u_{n}\|^{2}+\int_{\mathbb{R}^{4}}\frac{H(x, u_{n})}{|x|^{\eta}}dx  \\
             =&4I(u_{n})-I'(u_{n})u_{n}\\
             =&4I(u_{n})+o_{n}(1)
             =4c+o_{n}(1),\\
\end{aligned}
\right.
\end{equation*}
which is a contradiction to \eqref{(3.4)}. This proves that\ $\{u_{n}\}$ is bounded in $E$.    $\hfill\Box$\\
\noindent{\bf Proof of Theorem 1.2.} From Lemma 5.2, we have that the Cerami sequence\ $\{u_{n}\}$ is bounded in\ $E$. Applying  the same procedure in proof of Theorem 1.1, we will derive that $I'(u) = 0$ and $I(u) = c$. Moreover, we also get that\ $u$ is nonzero and $u$ is ground state. ~~~~~~~~$\hfill\Box$

\def\refname{References }

\end{document}